\documentclass{tac}
\usepackage{amsfonts,amssymb,mathtools}

\usepackage{tikz}
\usetikzlibrary{cd}

\newcommand{\op}{\textrm{op}}

\newcommand{\maps}{\colon}

\newcommand{\To}{\Rightarrow}
\newcommand{\simrightarrow}{\xrightarrow{\raisebox{-3pt}[0pt][0pt]{\ensuremath{\sim}}}}

\newcommand{\A}{\mathbb{A}}
\newcommand{\C}{\mathbb{C}}
\newcommand{\F}{\mathbb{F}}
\newcommand{\N}{\mathbb{N}}
\renewcommand{\P}{\mathbb{P}}

\newcommand{\Z}{\mathbb{Z}}

\newcommand{\Aut}{\mathrm{Aut}}

\newcommand{\Fin}{\mathsf{Fin}}
\newcommand{\Set}{\mathsf{Set}}
\newcommand{\Comm}{\mathsf{Comm}}
\newcommand{\Ring}{\mathsf{Ring}}

\newcommand{\Spec}{\mathsf{Sp}}
\newcommand{\core}{\mathsf{core}}

\newcommand{\CC}{\mathsf{C}}

\newcommand{\define}[1]{{\bf \boldmath{#1}}}

\title{Dirichlet Species and Arithmetic Zeta Functions}

\author{John\ C.\ Baez}

\address{School of Mathematics, University of Edinburgh, James Clerk Maxwell Building, Peter Guthrie Tait Road, Edinburgh, UK EH9 3FD.}

\eaddress{baez@math.ucr.edu}

\keywords{Day convolution, Dirichlet series, species, zeta function}
\amsclass{11M38,11R42,14G10,18M80}

\begin{document}

\maketitle

\begin{abstract}
Though Joyal's species are known to categorify generating functions in enumerative combinatorics, they also categorify zeta functions in algebraic geometry.  The reason is that any scheme $X$ of finite type over the integers gives a `zeta species' $Z_X$, and any species $F$ gives a Dirichlet series $\widehat{F}$, in such a way that $\widehat{Z}_X$ is the arithmetic zeta function of $X$, a well-known Dirichlet series that encodes the number of points of $X$ over each finite field.   Specifically, a $Z_X$-structure on a finite set is a way of making that set into a semisimple commutative ring, say $k$, and then choosing a $k$-point of the scheme $X$.  This is an elaboration of joint work with James Dolan.
\end{abstract}

\section{Introduction}
\label{sec:intro}

In combinatorics, structures on finite sets are often counted with the help of formal power series called generating functions \cite{FS,Wilf}.   Andr\'e Joyal categorified the theory of generating functions and showed that `species' could be used, not merely to count structures on finite sets, but also to work with them directly \cite{Joyal1,Joyal2}.  The business of counting solutions to polynomial equations over finite fields is often considered a separate topic from enumerative combinatorics, but species have a role to play here too.   In algebraic geometry one counts such solutions using  `arithmetic zeta functions',  which are not formal power series, but rather Dirichlet series \cite{Serre}.   However, we show here that any species gives a Dirichlet series, and all arithmetic zeta functions naturally arise from species in this way.   

For an example of an arithmetic zeta function, consider any collection of polynomial equations with integer coefficients in some finite set of variables.    For any commutative ring $k$, let $X(k)$ be the set of all solutions of these equations where the variables take values in $k$.  Any ring homomorphism $k \to k'$ determines a map $X(k) \to X(k')$, and it is easy to see this defines a functor $X \maps \Comm\Ring \to \Set$.    Any finite field has cardinality $p^n$ for some prime $p$, and  up to isomorphism there is a unique field of cardinality $p^n$, called $\F_{p^n}$.   The arithmetic zeta function of $X$ is defined to be
\begin{equation}
\label{eq:zeta_of_equations}
  \zeta_X(s) =  \prod_{p \in \P}  \exp \left( \sum_{n \ge 1} \frac{|X(\F_{p^n})|}{n} p^{-n s} \right)  
\end{equation}
where $\P$ is the set of prime numbers.   Thus, $\zeta_X$ encodes how many solutions the equations have in each finite field.

While Equation \eqref{eq:zeta_of_equations} looks curious at first sight, arithmetic zeta functions are central to many problems in number theory, from the Weil Conjectures proved by Grothendieck and Deligne \cite{Milne}, to the still unproved Birch and Swinnerton-Dyer conjectures \cite{Wiles}, and beyond.  Even the Riemann Hypothesis can be seen as a question about an arithmetic zeta function.  Indeed, consider the degenerate system of equations with just one solution in each commutative ring.   This gives the terminal object $T$ in the category of functors from $\Comm\Ring$ to $\Set$, and the arithmetic zeta function of $T$ is the Riemann zeta function $\zeta$, since
\[  \zeta_T(s) =  \prod_{p \in \P}  \exp \left( \sum_{n \ge 1} \frac{1}{n} p^{-n s} \right)  = \prod_{p \in \P} \frac{1}{1 - p^{-s}} = 
\sum_{n \ge 1} n^{-s} = \zeta(s) .\]

How can we obtain an arithmetic zeta function from a species?   Intuitively, a species is a structure one can put on finite sets and transport along bijections: for example, a coloring, or ordering, or the structure of being the vertices of a graph.  Suppose $F$ is some such structure.  If we denote the set of $F$-structures on a chosen $n$-element set by $F(n)$, the exponential generating function of $F$ is defined to be
\[    |F|(x) = \sum_{n \ge 0} \frac{|F(n)|}{n!} x^n. \]
This is a well-defined formal power series if  $F(n)$ is finite for all $n$.   In this case we can also associate a Dirichlet series to $F$, namely 
\[     \widehat{F}(s) = \sum_{n \ge 1} \frac{|F(n)|}{n!} n^{-s} . \]
This is again a formal series: we do not demand that it converges.

Here we show that arithmetic zeta functions are in fact the Dirichlet series of certain species.  Any functor
\[  X \maps \Comm\Ring \to \Set \]
gives a species $Z_X$ for which a $Z_X$-structure on a finite set is a way of making it into a semisimple commutative ring, say $k$, and then picking an element of $X(k)$.   In our main result, Theorem \ref{thm:zeta_functions}, we prove that when $X$ preserves products and maps finite fields to finite sets, we have
\[  \widehat{Z}_X(s) = \sum_{n \ge 1} \frac{|Z_X(n)|}{n!} n^{-s} \, = \,  \prod_{p \in \P}  \exp \left( \sum_{n \ge 1} \frac{|X(\F_{p^n})|}{n} p^{-n s} \right)  \]
where $\P$ is the set of prime numbers.    This clarifies the initially \emph{ad hoc} appearance of Equation \eqref{eq:zeta_of_equations}.  

Furthermore, in Theorem \ref{thm:zeta_species_3} we prove that the above equation between Dirichlet series arises from an isomorphism between species:
\[       Z_X \cong \prod_{p \in \P}^D  \exp_D F_{X,p}.  \]
where:
\begin{itemize}
\item $F_{X,p}$ is the species defined in Lemma \ref{lem:zeta_species_2}, for which an $F_{X,p}$-structure on a finite set is a way of making that set into a field of characteristic $p$, say $k$, and then choosing an element of $X(k)$.
\item $\exp_D$ is the Dirichlet exponential of a species, defined in Section \ref{subsec:Dirichlet_exponentiation}.
\item $\prod_{p \in P}^D$ is the Euler product of species, defined in Section \ref{subsec:Euler_products}.
\end{itemize}
Thus, our work categorifies a bit of the theory of arithmetic zeta functions, by showing that some basic equations arise from isomorphisms between species.

\section{The Riemann zeta function}
\label{sec:Riemann_zeta}

Before plunging into the general theory, let us consider the basic example: the Riemann zeta function.   We tackle this in a brute-force way and make only slight progress.  Everything becomes easier after we have introduced more technology.   Still, this first attempt is amusing and perhaps instructive.

Define a `$Z$-structure' on a finite set to be a way of giving it addition and multiplication operations that make it into a semisimple commutative ring.   We can count the $Z$-structures on an $n$-element set for a few small values of $n$.   To do this, we start by recalling two facts:

\begin{itemize}
\item By the Artin--Wedderburn theorem, a finite semisimple commutative ring is the same as a finite product of finite fields \cite[Thm.\ 5.2.4]{Cohn}.
\item Up to isomorphism there is one field with $q$ elements, denoted $\F_q$, when $q$ is a positive integer power of a prime number, and none otherwise \cite[Thm.\ 7.8.2]{Cohn}.
\end{itemize}

This lets us classify the semisimple commutative rings with $n$ elements `by hand' for small 
values of $n$.   For example:

\begin{itemize}
\item There are none when $n = 0$.
 
\item There is one when $n = 1$: the ring with one element.  (This is the \emph{empty} product of finite fields.)

\item There is one when $n = 2$: $\F_2$.

\item There is one when $n = 3$: $\F_3$.

\item There are two when $n = 4$: $\F_4$ and $\F_2 \times \F_2$.

\item There is one when $n = 5$: $\F_5$.

\item There is one when $n = 6$: $\F_2 \times \F_3$.

\item There is one when $n = 7$: $\F_7$.

\item There are three when $n = 8$: $\F_8$, $\F_2 \times \F_4$ and $\F_2 \times \F_2 \times \F_2$.
\end{itemize}

Next, how many ways are there to make an $n$-element set into a ring isomorphic to one on our list?  For starters: how many ways are there to make an $n$-element set into a ring isomorphic to some fixed $n$-element ring, say $k$?  Each bijection between the set and the ring $k$ gives a way to do this, but not all of them give different ways: two bijections give the same ring structure if and only if they differ by an automorphism of $k$.  So, the answer is $n!/|\Aut(k)|$.

To go further, we need to understand the automorphisms of finite fields \cite[Thm.\ 7.8.3]{Cohn}: if $q = p^m$ for some prime $p$, then 
\[  \Aut(\F_q) \cong \Z/m , \]
namely the cyclic group generated by the Frobenius automorphism 
\[  \begin{array}{cccl}
F \maps  &\F_q &\to& \F_q  \\
     &           x &\mapsto & x^p    \, .
\end{array}
\]
More generally, given a finite product of finite fields, its automorphisms all come from automorphisms of the factors together with permutations of like factors.  So, for example, $\F_2 \times \F_4$ has 2 automorphisms (coming from automorphisms of the second factor), while $\F_2 \times \F_2 \times \F_2$ has 6 (coming from permutations of the factors), and $\F_2 \times \F_2 \times \F_4$ has 4 (coming from permuting the first two factors, and automorphisms of the third).

Using these facts, we can count how many different ways there are to give an $n$-element set the structure of a semisimple commutative ring:
\begin{itemize}
\item For $n = 0$ there are $0$ ways.
 
\item For $n = 1$ there is $1!/1 = 1!$ ways.

\item For $n = 2$ there are $2!/1 = 2!$ ways.

\item For $n = 3$ there are $3!/1 = 3!$ ways.

\item For $n = 4$ there are $4!/2 + 4!/2 = 4!$ ways.

\item For $n = 5$ there are $5!/1 = 5!$ ways.

\item For $n = 6$ there are $6!/1 = 6!$ ways.

\item For $n = 7$ there are $7!/1 = 7!$ ways.

\item For $n = 8$ there are $8!/3 + 8!/2 + 8!/6 = 8!$ ways.
\end{itemize}
From this evidence, we might boldly guess there are always $n!$ ways, except for $n = 0$ when there
are none.  Though it is hard to see why from our work so far, we prove in Theorem \ref{thm:Riemann_zeta}
that this guess is correct.  In other words, the Dirichlet function of the Riemann species is
\[  \sum_{n \ge 1} \frac{n!}{n!} n^{-s} = \sum_{n \ge 1} n^{-s}  , \]
which is the Riemann zeta function $\zeta(s)$.

\section{Dirichlet species}

\subsection{The Dirichlet product of species}
\label{subsec:Dirichlet_product}

We call the groupoid of finite sets and bijections $\core(\Fin\Set)$, since it is the groupoid core of the category $\Fin\Set$, whose objects are finite sets and whose morphisms are functions.  The category of \define{species}  is the functor category $[\core(\Fin\Set), \Set ]$.  Any species, say $F \maps \core(\Fin\Set) \to \Set$, describes a type of structure that you can put on a finite set.  For any finite set $n$, $F(n)$ is the collection of structures of that type that you can put on the set $n$.

Here we are mainly interested in species where there are only finitely many structures of that type on any finite set.  So, we define the category of \define{tame species} to be the functor category $[\core(\Fin\Set), \Fin\Set]$.
Like the category of species itself, the category of tame species has quite a few interesting  symmetric monoidal structures.  Two of these come from coproducts and products in the target category $\Fin\Set$:

\begin{itemize}
\item the coproduct of species, given by 
\[ (F + G)(n) = F(n) + G(n) \]
This is usually called \define{addition}.
\item the product of species, given by
\[  (F \times G)(n) = F(n) \times G(n) \]
This is often called the \define{Hadamard product}.
\end{itemize}
Two more symmetric monoidal structures arise via Day convolution.  Recall that any symmetric monoidal structure on a category $\CC$ extends uniquely one on $[\CC^\op, \Set]$ such that the tensor product preserves colimits in each argument \cite{Day}.  Since this new monoidal structure is analogous to convolving functions on a group, it is called `Day convolution'.   While $\core(\Fin\Set)^\op \cong\core(\Fin\Set)$ does not have coproducts or products, it inherits symmetric monoidal structures which we call $+$ and $\times$ from $\Fin\Set$, which does.   By Day convolution each of these monoidal structures on $\core(\Fin\Set)$ gives a symmetric monoidal structure on the category of species, which restricts to one on tame species:

\begin{itemize}
\item The symmetric monoidal structure on species arising via Day convolution from the $+$ monoidal
structure on $\core(\Fin\Set)$ is called the \define{Cauchy product} and denoted $\cdot_C$.
\item The symmetric monoidal structure on species arising via Day convolution from the $\times$ monoidal
structure on $\core(\Fin\Set)$ is called the \define{Dirichlet product} (or \define{arithmetic product}) and
denoted $\cdot_D$.
\end{itemize}

To understand the Dirichlet product it is useful to warm up with the more familiar Cauchy product.
The rough idea is that if $F$ and $G$ are two species, an $F \cdot_C G$-structure on a finite set $S$ is a way of writing $S$ as a coproduct of two sets and putting an $F$-structure on the first set and a $G$-structure on the second.   We can make this precise and avoid `overcounting' the coproduct decompositions as follows:
\begin{equation}
\label{eq:Cauchy_product}
 (F \cdot_{C} G)(S) = \sum_{T \subseteq S} F(T) \times G(S-T) 
 \end{equation}
where $S$ is a finite set and we take the coproduct over all subsets $T \subseteq S$.   
The Dirichlet product is in some sense dual.   Now the idea is that an $F \cdot_D G$-structure on a finite set is a way of writing $S$ as a product of two sets and putting an $F$-structure on the first set and a $G$-structure on the second.    Maia and M\'endez \cite{MM} gave a convenient way to make this precise.  They defined a \define{cartesian decomposition} of a set $S$ to be an ordered pair $(\pi_1, \pi_2)$ of equivalence relations on $S$ such that each equivalence class of $\pi_1$ intersects each equivalence class of $\pi_2$ in a singleton.   If we let $S_i$ be the set of equivalence classes of $\pi_i$, then the maps sending elements of $S$ to their equivalence classes give a product cone:
\[
\begin{tikzcd}[column sep = small]
    & S \arrow[dl] \arrow[dr]
    \\
S_1  & & S_2.
\end{tikzcd}
\] 
The  Dirichlet product of species $F$ and $G$ is then given by
\begin{equation}
\label{eq:Dirichlet_product}
(F \cdot_D G)(S) = \sum_{\text{cartesian decompositions } (\pi_1,\pi_2) \text{ of } S} F(S_1) \times G(S_2) .
\end{equation}
A bijection $f \maps S \to S'$ induces a bijection between cartesian decompositions of $S$ and cartesian decompositions of $S'$, and thus between $(F \cdot_D G)(S)$ and $(F \cdot_D G)(S')$.  Using this, $F \cdot_D G$ becomes a species.

We can also visualize a cartesian decomposition of a finite set as equivalence class of ways of organizing its elements into a rectangle, where two ways are equivalent if they differ by permuting rows and/or columns, for example
\[  
\begin{array}{ccc}
1 & 4 & 3 \\
2 & 6 & 5 
\end{array} 
\; \sim \;
\begin{array}{ccl}
6 & 5 & 2 \\
4 & 3 & 1. 
\end{array} 
\]
Dwyer and Hess \cite{DH} rediscovered the Dirichlet product by considering the monoidal structure on `symmetric sequences', which form a category equivalent to that of species.   For more, see the work of Gambino, Garner and Vasilakopoulou \cite{GGV}.

There are, of course,
\[          \binom{n}{k} = \frac{n!}{k! (n - k)!}. \]
$k$-element subsets of the ordinal $n$.   This makes it evident that when we take cardinalities, the Cauchy product obeys
\[   |(F \cdot_{C} G)(n)| = \sum_{k = 0}^n \, \binom{n}{k}  \, |F(k)|\, |G(n - k)| .\]
The Dirichlet product obeys a similar formula.   A cartesian decomposition of an $n$-element set into a $k$-element set and some other set exists if and only $k$ divides $n$.  In this case the group $S_n$ acts transitively on the set of such cartesian decompositions, and the stabilizer of any one is 
$S_k \times S_{n/k}$, so the number of such cartesian decompositions is
 \[    {n \brace k} = \frac{n!}{k! (n/k)!}. \]
The notation here was introduced by  Maia and M\'endez \cite[Sec.\ 2]{MM}.   When we take cardinalities, we thus have
\begin{equation}
\label{eq:dirichlet_product_1}
  |(F \cdot_{D} G)(n)| = \sum_{k | n} {n \brace k} \, |F(k)|\, |G(n/k)| . 
\end{equation}
where we sum over natural numbers $k$ that divide $n$.    
\subsection{The category of Dirichlet species}
\label{subsec:dirichlet_species}

The sum in Equation \eqref{eq:dirichlet_product_1} 
is finite except when $n = 0$, since every natural number divides $0$.  This wrinkle
means that the Dirichlet product of tame species may not be tame!   However, the Dirichlet product of 
two tame species that vanish on the empty set is tame, and it vanishes on the empty set.  This leads us to the following definition:

\begin{definition}  Define the category of species to be 
\[     \Spec = [\core(\Fin\Set), \Set] \]
and define the category of \define{Dirichlet species}, $\Spec_D$, to be the full subcategory of $\Spec$ whose objects are species $F$ with $F(\emptyset) = \emptyset$.
\end{definition}  

Equivalently, we can think of a Dirichlet species as a functor $F \maps \core(\Fin\Set_{> 0}) \to \Set$
where $\Fin\Set_{> 0}$ is the category of \emph{nonempty} finite sets and functions between them.
Indeed, we have
\[   \Spec_D \simeq [\core(\Fin\Set_{> 0}), \Set]  \]
and this is a useful way to think about the category of Dirichlet species.

To spell this out a bit more, note that the inclusion
\[        i \maps \core(\Fin\Set_{> 0}) \to \core(\Fin\Set)   \]
is symmetric monoidal with respect to the $\times$ monoidal structure, so by Day convolution it
gives a symmetric monoidal functor $i^\ast \maps \Spec \to \Spec_D$ which takes any species and restricts it to nonempty sets.

\begin{proposition}
The symmetric monoidal functor $i^\ast  \maps \Spec \to \Spec_D$ has a left adjoint $i_! \maps \Spec_D \to \Spec$ which extends any Dirichlet species $F$ to a species $i_!F$ with $(i_!F)(\emptyset) = \emptyset$.
This functor $i_!$ is also symmetric monoidal, and the unit of the adjunction
\[   \eta \maps 1 \To  i^\ast \circ i_! \]
is a symmetric monoidal natural isomorphism, so $\Spec_D$ is not only a coreflective
subcategory of $\Spec$, but a `coreflective symmetric monoidal subcategory'.
\end{proposition}

\begin{proof}
Given a Dirichlet species $F$ and a species $G$, any morphism of species $i_! F\to G$ restricts to
a morphism $F \to i^\ast G$ of Dirichlet species.  On the other hand, any morphism $F \to i^*G$ of Dirichlet species  extends uniquely to a morphism $i_! F \to G$ of species, since $(i_!F)(\emptyset)$ is initial.  
These operations determine a bijection $\Spec(i_!F, G) \cong \Spec_D(F, i^*G)$, natural in $F$ and $G$, 
so $i_!$ is left adjoint to $i^\ast$.    The unit $\eta$  of this adjunction is the natural map from a Dirichlet species $F$ to the Dirichlet species formed by first extending $F$ to a species $i_!F$ and then restricting it back to the groupoid of nonempty finite sets; this is clearly an isomorphism.

While the left adjoint of a symmetric monoidal functor is in general only symmetric \emph{oplax} monoidal, 
it can be seen by direct computation that $i_!$ is strong monoidal, and that the counit and unit of the adjunction between $i^\ast$ and $i_!$ are symmetric monoidal natural transformations.
\end{proof}

 \subsection{Dirichlet series}
 \label{subsec:Dirichlet_series}
 
We want to define operations on species that categorify familiar operations on Dirichlet series.  To set these up we need to understand how a tame species gives a Dirichlet series.
 
\begin{definition}
\label{defn:Dirichlet_series}
Given a tame species $F$, we define its \define{Dirichlet series} by
\[ \widehat{F}(s) = \sum_{n \ge 1} \frac{|F(n)|}{n!} n^{-s}  \]
where the sum is a formal one.
\end{definition}

It is easy to check that if $F$ and $G$ are tame species, then
\[  \widehat{\phantom{\big|} F+G \phantom{\big|}} = \widehat{F} + \widehat{G}. \]
More interestingly, if $F$ and $G$ are tame Dirichlet species, then
\[ \widehat{\phantom{\big|} F \cdot_{D} G \phantom{\big|}} = \widehat{F} \, \widehat{G} .\]
This was shown by Maia and M\'endez \cite[Prop.\ 2.3]{MM}, but the proof is so quick we give it here:
\[   \begin{array}{ccl}
 \widehat{\phantom{\big|}F \cdot_{D} G \phantom{\big|}}(s) &=&
 \displaystyle{  \sum_{n \ge 1} \frac{|(F \cdot_{D} G)(n)|}{n!} n^{-s} } \\[5 pt]
 &=& \displaystyle{  \sum_{n \ge 1} \sum_{k | n} {n \brace k} \frac{|F(k)|\, |G(n/k)|}{n!} n^{-s}   }  \\[5pt]
  &=& \displaystyle{  \sum_{n \ge 1} \sum_{k | n} \frac{|F(k)|}{k!} k^{-s} \; \frac{|G(n/k)|}{(n/k)!} (n/k)^{-s}   }  \\[15pt]
  &=& \widehat{F}(s) \, \widehat{G}(s).
\end{array}
\]
where in the second step we use Equation \eqref{eq:dirichlet_product_1}.

All of this resembles a more familiar story involving generating functions and the Cauchy product.  Namely, any tame species $F$ has an \define{exponential generating function}
\[ |F|(z) = \sum_{n \ge 0} \frac{|F(n)|}{n!} z^n .\]
It is well known \cite{BLL,Joyal1} that
\[  |F+G| = |F| + |G| \]
and 
\[ |F \cdot_{C} G| = |F|\, |G|. \]
The latter follows from a calculation very much like the one given above for the Dirichlet product, but with binomial coefficients replacing Dirichlet binomials.

In short, Dirichlet series are adapted to the Dirichlet product of species just as exponential generating functions are adapted to the Cauchy product.   For Dirichlet species they present the same information, since the Dirichlet series $\widehat{F}$ is obtained from the generating function $|F|$ by the replacement
\[ z^n \mapsto \left\{ \begin{array}{ccl} n^{-s} & n \ge 1 \\ [3pt]
                                                                      0 & n = 0 . 
                                                                       \end{array}  \right.\]
In number theory a closely related map is called the Mellin transform, and this explains the frequent appearance of Mellin transforms in number theory.  

\subsection{Dirichlet exponentiation}
\label{subsec:Dirichlet_exponentiation}

Any Dirichlet series 
\[                    f(s) = \sum_{n \ge 1} \frac{f(n)}{n!} n^{-s}   \]
has an exponential
\[                   \exp(f)(s) = \sum_{m \ge 0} \frac{1}{m!} \left(  \sum_{n \ge 1} \frac{f(n)}{n!} n^{-s}    \right)^m \]
which is again a well-defined Dirichlet series.   Suppose $F$ is a tame Dirichlet species.  We would like 
to define a tame species $\exp_D(F)$ such that 
\[    \widehat{\phantom{\large|} \exp_D(F) \phantom{\large|}} = \exp(\widehat{F}) .\]
However, this is impossible unless $F$ vanishes on one-element sets.   
To get a feeling for this, consider the unit object for the Dirichlet product, $I$.   This 
Dirichlet species deserves to be called \define{being a one-element set}, since it is defined by
\begin{equation}
\label{eq:unit_object}
 I(S) = \left\{ \begin{array}{ccl} 1 & |S| = 1 \\ 
                                                            \emptyset & \text{otherwise}  
                                                                       \end{array}  \right. 
\end{equation}                                                                   
where $1$ here stands for the ordinal, which is both a natural number and a one-element set.
We have $\widehat{I}(s) = 1$ and thus
\[   \exp(\widehat{I})(s)= e. \]
This cannot be the Dirichlet series of a tame species, since the only constants that are Dirichlet series of
tame species are natural numbers.

We can, however, find  a tame species $\exp_D(F)$ whose Dirichlet series is $\exp(\widehat{F})$ when $F$ vanishes on one-element sets.  The rough idea is that a $\exp_D(F)$-structure on a finite set is a way of writing that set as an unordered product of finite sets and putting an $F$-structure on each factor.   

One way to make this idea precise is to consider, for each $n \in \N$, the $n$-fold Dirichlet product of $F$ with itself, which we call the $n$th \define{Dirichlet power} of $F$:
\[    F^n_D = \underbrace{F \cdot_D F \cdot_D \cdots \cdot_D F}_{n \text{ factors }}. \]
(We define $F^0_D$ to be $I$, the unit object for the Dirichlet product.)
Since the Dirichlet product is symmetric monoidal, the symmetric group $S_n$ acts on $F^n_D$.  Since the category of Dirichlet species is a presheaf category, it has colimits, so we can take the quotient of $F^n_D$ by this $S_n$ action and get a Dirichlet species that we call the \define{symmetrized $n$th Dirichlet power} and denote as
\[   F_D^n/S_n .\]
Finally, we make the following definition:

\begin{definition}
\label{defn:Dirichlet_exponential}
Given a Dirichlet species $F$, we define its \define{Dirichlet exponential} to be the coproduct
\[    \exp_D(F) = \sum_{n \ge 0}  F_D^n/S_n \]
in the category of Dirichlet species.
\end{definition}

\begin{lemma}
\label{lem:Dirichlet_exponential}
For any tame Dirichlet species $F$ for which $F(1) = \emptyset$, the species $\exp_D(F)$ is tame and
\[    \widehat{\phantom{\large|} \exp_D(F) \phantom{\large|}} = \exp(\widehat{F}) .\]
\end{lemma}

\begin{proof}
Since for any two Dirichlet species $F$ and $G$ we have 
\[ \widehat{\phantom{\big|} F \cdot_{D} G \phantom{\big|}} = \widehat{F} \, \widehat{G} ,\]
it follows by induction that the Dirichlet powers obey
\begin{equation}
 \label{eq:Dirichlet_powers}
  \widehat{F_D^n} = (\widehat{F})^n .
 \end{equation}
Since colimits in a presheaf category are computed pointwise, for any finite set $T$ we have
\[        (F_D^n/S_n) (T) = (F_D^n(T))/S_n  .\]
We can understand the action of $S_n$ on the set $F_D^n(T)$ using the description of the Dirichlet product in terms of cartesian decompositions in Section \ref{subsec:Dirichlet_product}.   An element of $F_D^n(T)$ is a way of writing $T$ as a cartesian product of $n$ sets and choosing an $F$-structure on each of these sets.   More precisely, it is a way of 
\begin{itemize}
\item choosing $n$ equivalence relations $\pi_1 , \dots, \pi_n$ on $T$ such that if $T_i \subseteq T$ is any equivalence class for $\pi_i$ then $T_1 \cap \cdots \cap T_n$ is a singleton;
\item choosing an $F$-structure $x_i$ on the set of equivalence classes of $\pi_i$ for each $1 \le i \le n$.
\end{itemize}
The group $S_n$ acts on $F_D^n(T)$ by permuting the equivalence relations $\pi_1, \dots, \pi_n$ along with
the corresponding $F$-structures $x_1, \dots, x_n$.   This action is free if for each $i$ the set of equivalence classes of $\pi_i$ has more than one element, since in that case \(i \ne j\) implies \(\pi_i \ne \pi_j\).   Thus, since $F$ is empty for one-element sets, we have
\begin{equation}
 \label{eq:symmetrized_Dirichlet_powers}
    |(F_D^n/S_n)(T)| = \frac{1}{n!} |F_D^n(T)|   .
 \end{equation}
We thus have
\[  \begin{array}{ccll}
 \widehat{\phantom{\big|} \exp_D(F) \phantom{\big|}}(s) &=& \displaystyle{ \sum_{n \ge 0} \widehat{\phantom{\big|}F_D^n/S_n \phantom{\big|}}(s) } &
 \text{by Definition \ref{defn:Dirichlet_exponential} }  \\ \\
 &=& \displaystyle{ \sum_{n \ge 0}  \sum_{m \ge 1} \frac{|(F_D^n/S_n)(m)|}{m!} m^{-s} }& 
 \text{by Definition \ref{defn:Dirichlet_series} }\\ \\ 
 &=& \displaystyle{ \sum_{n \ge 0}  \sum_{m \ge 1} \frac{|F_D^n(m)|}{n!m!} m^{-s}} 
  & \text{by Equation (\ref{eq:symmetrized_Dirichlet_powers}) }\\ \\
  &=& \displaystyle{ \sum_{n \ge 0} \frac{1}{n!} \widehat{F_D^n}(s)} &
   \text{by Definition \ref{defn:Dirichlet_series} }  \\ \\
&=& \displaystyle{ \sum_{n \ge 0} \frac{1}{n!} \widehat{F}(s)^n}  & 
\text{by Equation (\ref{eq:Dirichlet_powers}) }
\\ \\
&=& \exp(\widehat{F})(s) .  
 \end{array}
 \]
\end{proof}

Lemma \ref{lem:Dirichlet_exponential} should remind us of a similar result for the `usual' exponential of a species $F$, which we should perhaps call the \textbf{Cauchy exponential} $\exp_C(F)$. This is a species such that an $\exp_C(F)$-structure on a finite set is the same as a way of writing that set as an unordered coproduct of finite sets and putting an $F$-structure on each summand.   When $F$ is tame and $F(\emptyset) = \emptyset$, the Cauchy exponential is tame and
\[    |\exp_C(F)| = \exp(|F|) \, . \]
When $F(\emptyset)$ is nonempty, the Cauchy exponential is not tame and this equation fails, for reasons
analogous to the problem we saw for the Dirichlet exponential when $F(1)$ is nonempty.  The reason is that $\emptyset$ is the unit for the coproduct of finite sets, so we can write any finite set as a coproduct involving arbitrarily many copies of $\emptyset$, just as $1$ is the unit for the product, so we can write any finite set as a product involving arbitrarily many copies of $1$.

Here is an example important for understanding zeta functions:

\begin{example}
\label{ex:F_p}
Suppose $F_p$ is the species such that $F_p(S)$ is the set of ways of making the finite set $S$ into a field of characteristic $p$.   There are only finitely many ways to do this, and none unless $|S| = p^n$ for $n \ge 1$, so $F_p$ is a tame Dirichlet species.   When $n \ge 1$ there are $p^n!/n$ ways to make a $p^n$-element set into a field of characteristic $p$, since there are $p^n!$ bijections between the set and $\F_{p^n}$, and this field has $n$ automorphisms.  Thus, we have
\[ 
\widehat{F_p}(s) = \displaystyle{ \sum_{j \ge 1} \frac{|F_p(j)|}{j!} j^{-s} } 
=   \displaystyle{  \sum_{n \ge 1} \frac{|F_p(p^n)|}{p^n!} p^{-ns} }
=   \displaystyle{  \sum_{n \ge 1}\frac{1}{n} p^{-ns} } 
= \ln\left(\frac{1}{1 - p^{-s}}\right).
\]

The Dirichlet exponential $\exp_D(F_p)$ is the species such that $\exp_D(F_p)(S)$ is the set of ways of making the finite set $S$ into a product of fields of characteristic $p$.   
By Lemma \ref{lem:Dirichlet_exponential}, the Dirichlet series of this species is
\[     \widehat{\phantom{\big|}\exp_D(F_p) \phantom{\big|}}(s) 
 = \exp(\widehat{F_p}(s)) = \frac{1}{1 - p^{-s}}.\]
\end{example}

\subsection{Euler products}
\label{subsec:Euler_products}

In number theory, an \define{arithmetic function} is a function
\[ f \maps \N_{>0} \to \C \]
where $\N_{> 0} = \{1, 2, 3, \dots\}$.
An arithmetic function is said to be \define{multiplicative} if $f(1) = 1$ and $f(m n) = f(m) f(n)$ whenever $m$ and $n$ are relatively prime.  A multiplicative function is determined by its values on prime powers, so we may think of it as a function of isomorphism classes of finite fields.

Any arithmetic function $f$ gives a Dirichlet series 
\[  F(s) = \sum_{n \ge 1} f(n) n^{-s} \]
and $f$ is multiplicative if and only if this Dirichlet series has an \define{Euler factorization}:
\[  F(s) = \prod_{p \in \P} (1 + f(p) p^{-s} + f(p^2) p^{-2s} + \cdots ) \]
where $\P$ is the set of all primes.   Each factor here is a Dirichlet series in its own right, and 
the $p$th factor is `$p$-local':

\begin{definition}
A Dirichlet series $\sum_{n \ge 1} f(n) n^{-s}$ is \define{$p$-local} for some prime $p$ if $f(n) = 0$ unless $n$ is a power of $p$ (possibly the zeroth power), and $f(1) = 1$.   
\end{definition}

Conversely, suppose for each prime $p$ we have a $p$-local Dirichlet series $F_p$.   Then we can define the infinite product
\[      F(s) = \prod_{p \in \P} F_p(s)   \]
as a Dirichlet series, and if we write
\[  F(s) = \sum_{n \ge 1} f(n) n^{-s} \]
then $f$ is a multiplicative arithmetic function.

All this is ripe for categorification.  

\begin{definition} 
A Dirichlet species $A$ is \define{$p$-local} for some prime $p$ if $\widehat{A}$ is $p$-local: that is,
$A(S)$ is empty unless $|S|$ is a power of $p$, and $A(S)$ has one element if $|S| = 1$.
\end{definition}

Suppose for each prime $p$ we have a $p$-local Dirichlet species $A_p$.    Then we can take the Dirichlet product of all the $A_p$ and get a well-defined Dirichlet species.  Since this is an infinite product, its formulation requires some care.  Rather than developing a general theory of infinite Dirichlet products, it seems quicker to leverage what we know about finite Dirichlet products.

Let $p_i$ be the $i$th prime number.   Let $A[k]$ be the Dirichlet product of the species $A_p$ for the first $k$ primes:
\[      A[k] = A_2 \cdot_D A_3 \cdot_D A_5 \cdot_D \cdots \cdot_D A_{p_k} .\]
For any Dirichlet species $X$ there is a unique morphism of Dirichlet species from $I$ to $X$ where $I$ is the unit for the Dirichlet product described in Equation \eqref{eq:unit_object}.    Tensoring this with the identity morphism on $A[k]$ gives a morphism 
\[         \iota_k \maps A[k] \to A[k+1].   \]
We use these morphisms $\iota_k$ to define the \define{Euler product} of the Dirichlet species $A_p$ as follows:
\[ {\prod_{p \in \P}}^D \,  A_p   = \underset{k \to \infty}{\mathrm{colim}} \, A[k]   .\]
Since the category of Dirichlet species is cocomplete, this is a well-defined Dirichlet species.

\begin{lemma}  
\label{lem:Euler_product}
If $A_p$ are tame $p$-local Dirichlet species, one for each $p \in \P$, then their Euler product
\[ A = {\prod_{p \in \P}}^D \,  A_p  \]
is a tame Dirichlet species, with
\[    \widehat{A} = \prod_{p \in \P} \widehat{A}_p ,\]
and the arithmetic function $n \mapsto |A(n)|$ is multiplicative.
\end{lemma}

\begin{proof}
Suppose $S$ is a finite set.   By definition, an element of $A[k](S)$ consists of:
\begin{itemize}
\item a list of equivalence relations $\pi_1, \dots, \pi_k$ on $S$ such that if $S_i$ is the set of equivalence classes of $\pi_i$, the quotient maps $\pi_i \maps S \to S_i$ form a product cone expressing $S$ as the product $S_1 \times \cdots \times S_k$,
\end{itemize}
together with
\begin{itemize}
\item an element $x_i \in A_{p_i}(S_i)$ for each $i = 1, \dots, k$.
\end{itemize}
Since each species $A_p$ is $p$-local, $A_{p_i}(S_i) = \emptyset$ unless $|S_i|$ is a power of $p_i$, and $|A_{p_i}(S_i)| = 1$ if $|S_i| = 1$. Thus $A[k](S) = \emptyset$ whenever any prime factor of $|S|$ is not among the primes $p_1, \dots, p_k$, but whenever $k$ is large enough that all the prime factors of $|S|$  are among these primes, the map
\[    \iota_k(S) \maps A[k](S) \to A[k+1](S)  \]
is a bijection.

It follows that 
\begin{equation}
\label{eq:isomorphism}
A(S) \cong A[k](S) 
\end{equation}
whenever $k$ is large enough that all the prime factors of $|S|$ are among the first $k$ primes.
Since the sets $A_p(S)$ are finite for each prime $p$, so is each set $A[k](S)$, and thus so is $A[S]$.  Thus
$A$ is tame.

By Equation \eqref{eq:isomorphism}, the first $n$ terms of the Dirichlet series $\widehat{A}$ and $\widehat{A}[k]$ are equal when $n \le p_k$.    Since 
\[   \widehat{A[k]}(s) = \prod_{i = 1}^k \widehat{A_{p_i}}(s)   \]
it follows that 
\[   \widehat{A}(s) = \prod_{p \in \P} \widehat{A_p}(s) . \]
Furthermore, since the arithmetic function $n \mapsto |A[k](n)|$ is multiplicative for any $k$, 
and $|A(n)| = |A[k](n)|$ when $n \le p_k$, it follows that $n \mapsto |A(n)|$ is multiplicative as well.  \end{proof}

We can now finally prove our claim in Section \ref{sec:Riemann_zeta}, namely that the Riemann zeta function is the Dirichlet series of the species of semisimple commutative rings.  

Recall from Example \ref{ex:F_p} that an $\exp_D(F_p)$-structure on a finite set $S$ is a way of making it into a product of fields of characteristic $p$.   There is no way to do this unless $|S|$ is a power of $p$, and one way if $|S| = 1$, since a one-element ring is an empty product of fields.  Thus $\exp_D(F_p)$ is a $p$-local Dirichlet species.  Since each Dirichlet species $F_p$ is tame, so is $\exp_D(F_p)$.  Any finite semisimple commutative ring $k$ can be uniquely expressed as a product of rings $k_p$, one for each $p \in \P$, that are products of finite fields of characteristic $p$.   Thus if we define the \define{Riemann species} to be
\begin{equation}
\label{eq:Riemann_species}
  Z = {\prod_{p \in \P}}^D  \exp_D(F_p),  
\end{equation}
it follows that a $Z$-structure on a finite set is a way of making it into a semisimple commutative ring.   By Lemma \ref{lem:Euler_product}, $Z$ is a tame Dirichlet species.

We can now show that the Riemann species $Z$ has the Riemann zeta function $\zeta$ as its Dirichlet series.

\begin{theorem}
\label{thm:Riemann_zeta}
$\widehat{Z} = \zeta$.
\end{theorem}

\begin{proof}
In Example \ref{ex:F_p} we saw 
\[     \widehat{\phantom{\big|}\exp_D(F_p) \phantom{\big|}}(s) = \frac{1}{1 - p^{-s}} \]
so by Lemma \ref{lem:Euler_product} we have
\[    \widehat{Z}(s) = \prod_{p \in \P}  \widehat{\phantom{\big|}\exp_D(F_p) \phantom{\big|}}(s) = 
\prod_{p \in \P} \frac{1}{1 - p^{-s}} = \zeta(s) .\]
\end{proof}

\section{Zeta species}
\label{sec:zeta_species}

Algebraic geometry is often formulated in terms of schemes, but for the present purposes we can work with a simpler and more general concept, namely a functor
\[ X \maps \Comm\Ring \to \Set  . \]

For example:

\begin{itemize}
\item We can start with a polynomial equation with integer coefficients, say 
\[       x^3 + y^3 = z^3 \]
and let $X(k)$ be the set of solutions where the variables take values in the commutative ring $k$.  Any homomorphism $k \to k'$ gives a map $X(k) \to X(k')$, and it is easy to check that $X$ is a functor.  We call $X(k)$ the \define{set of $k$-points} of $X$.  
\item More generally, we can specify a functor $X \maps \Comm\Ring \to \Set$ by giving a commutative ring $R$ and letting $X(k) = \hom(R,k)$ be the set of ring homomorphisms $f \maps R \to k$.   The previous example is the special case where we take $R = \Z[x,y,z]/\langle x^3 + y^3 = z^3 \rangle$.  
\item More generally still, any scheme determines a functor from $\Comm\Ring$ to $\Set$, often called its `functor of points'.  The previous example is a special case of this: a so-called `affine scheme'.
\end{itemize}

\begin{definition} 
Given a functor $X \maps \Comm\Ring \to \Set$, we define its \define{zeta species}, $Z_X$, as follows: a $Z_X$-structure on a finite set is a way to make that set into a semisimple commutative ring, say $k$, and then choose an element of $X(k)$.  
\end{definition}

 Note also that the process of getting a zeta species from a functor $X \maps \Comm\Ring \to \Set$ is itself functorial.  That is, there is a functor
\[                Z \maps [\Comm\Ring, \Set] \to [\core(\Fin\Set), \Set]  \]
sending any functor $X \maps \Comm\Ring \to \Set$ to its zeta species $Z_X$, and any natural transformation between such functors to the evident natural transformation between their zeta species.

\begin{definition}  We say a functor $X \maps \Comm\Ring \to \Set$ is \define{tame} if $X(k)$ is finite whenever $k$ is a finite semisimple commutative ring. 
\end{definition}

Clearly a functor $X \maps \Comm\Ring \to \Set$ is tame if and only if its zeta species $Z_X$ is tame. 
Whenever $X$ is the functor of points of a scheme of finite type defined over the integers, it is tame \cite{Serre}. 

\begin{definition} A functor $X \maps \Comm\Ring \to \Set$ is \define{multiplicative} if it preserves products.
\end{definition}

Whenever $X$ is the functor of points of a scheme, it is multiplicative.  This is easiest to see for an affine scheme, since then $X = \hom(R, -)$ and any representable functor preserves products.  

\begin{lemma} 
\label{lem:zeta_species_1}
If $X \maps \Comm\Ring \to \Set$ is multiplicative, its zeta species $Z_X$ is an Euler product
\[                 Z_X \cong {\prod_{p \in \P}}^D \, Z_{X,p}  \]
where a $Z_{X,p}$-structure on a finite set is a way to make that set into a ring $k_p$ that is a product of finite fields of characteristic $p$ and then choose an element of $X(k_p)$. 
\end{lemma}

\begin{proof} 
To put a $Z_X$-structure on a finite set is to make it into a semisimple commutative ring $k$ and choose an element of $X(k)$.   But a finite semisimple commutative ring is the same as  a product over $p \in \P$ of rings $k_p$ where $k_p$ is a finite product of finite fields of characteristic $p$, and all but finitely many of the $k_p$ have just one element.  Since $X$ is multiplicative, to choose an element of $X(k)$ is the same as to choose an element of $X(k_p)$ for each $p$.  
\end{proof}

\begin{lemma}
 \label{lem:zeta_species_2}
 If $X \maps \Comm\Ring \to \Set$ is multiplicative, each species $Z_{X,p}$ can be written as the Dirichlet exponential
\[                Z_{X,p} \cong \exp_D \, (F_{X,p})  \]
where an $F_{X,p}$-structure on a finite set is a way to make that set into a field $k$ of characteristic $p$ and then choose an element of $X(k)$.
\end{lemma}

\begin{proof}
To put a $\exp_D (F_{X,p})$-structure on a finite set is to write it as a product of finite fields of characteristic $p$, say $k = f_1 \times \cdots \times f_n$, and then choose an element of $X(f_i)$ for each $i = 1, \dots, n$.    Since $X$ is multiplicative, this collection of choices is the same as choosing an element of $X(k)$.   Thus $\exp_D F(X_{X,p})$ is naturally isomorphic to $Z_{X,p}$.
\end{proof}

\begin{theorem}
\label{thm:zeta_species_3}
If $X \maps \Comm\Ring \to \Set$ is multiplicative, its zeta species is given by
\[                 Z_X \cong {\prod_{p \in \P}}^D   \exp_D \, (F_{X,p}) . \]
\end{theorem}

\begin{proof}
This follows straight from Lemmas \ref{lem:zeta_species_1} and  \ref{lem:zeta_species_2}.
\end{proof}

\section{Arithmetic Zeta Functions}
\label{sec:zeta_functions}

We are now ready to prove our main result, Theorem \ref{thm:zeta_functions}.  This says that we have a way of `categorifying' arithmetic zeta functions: that is, finding species with these zeta functions as their Dirichlet series.  We have seen that any functor $X \maps \Comm\Ring \to \Set$ gives a Dirichlet species $Z_X$ for which a $Z_X$-structure on a nonempty finite set is a way to make that set into a semisimple commutative ring, say $k$, and then choose an element of $X(k)$.   When $Z_X$ is tame, it has a Dirichlet series
\[    \widehat{Z}_X(s) = \sum_{n \ge 1} \frac{|Z_X(n)|}{n!} n^{-s}. \]
We now show that when $Z_X$ is also multiplicative, $\widehat{Z}_X$ equals the \define{arithmetic zeta function} of $X$, a Dirichlet series already familiar to algebraic geometers when $X$ is a scheme of finite type over the integers \cite{Serre}.  This is given by the right side of the equation below.

\begin{theorem} 
\label{thm:zeta_functions}
 Suppose a functor $X \maps \Comm\Ring \to \Set$ is tame and multiplicative.  Then 
\[    \widehat{Z}_X(s) =  \prod_{p \in \P}  \exp \left( \sum_{n \ge 1} \frac{|X(\F_{p^n})|}{n} p^{-n s} \right)  .\]
\end{theorem}

\begin{proof} 
By  Theorem \ref{thm:zeta_species_3} we have an isomorphism of tame Dirichlet species
\[                 Z_X \cong {\prod_{p \in \P}}^D   \exp_D \, (F_{X,p})  \]
so taking Dirichlet series and using Lemmas \ref{lem:Dirichlet_exponential} and \ref{lem:Euler_product}, we obtain
\[  \widehat{Z}_X(s) =  \prod_p \exp \left( \sum_{n \ge 1} \frac{|F_{X,p}(p^n)|}{p^n!}  p^{-n s} \right) \, . \]
It thus suffices to show that
\[    \frac{|F_{X,p}(p^n)|}{p^n!} =  \frac{|X(\F_{p^n})|}{n} \, . \]
By definition $|F_{X,p}(p^n)|$ is the number of ways to make an $p^n$-element set into a field, say $k$, and choose an element of $X(k)$.  But as we have seen, there are $p^n! / n$ ways to make a $p^n$-element set into a field.  So, there are $\frac{p^n!}{n} |X(\F_{p^n})| $ ways to make a $p^n$-element set into a field $k$ and then choose an element of $X(k) \cong X(\F_{p^n})$.  This gives the equation above.
\end{proof}

Again, the point is not so much the formula in Theorem \ref{thm:zeta_functions}: it is that Theorem \ref{thm:zeta_species_3} \emph{categorifies} this formula, lifting it to an isomorphism of species. 

\begin{example}
Let $T$ be the terminal object in $[\Comm\Ring,\Set]$: that is, the functor, unique up to isomorphism, such that $T(k)$ is a one-element set for every commutative ring $k$.   A $Z_T$-structure on a finite set is a way of making it into a semisimple commutative ring $k$, so $Z_T$ is isomorphic to the Riemann species $Z$ as defined in Equation \eqref{eq:Riemann_species}.   The functor $T$ is tame and multiplicative, so Theorem \ref{thm:zeta_functions} implies
\[      \widehat{Z}(s) = \widehat{Z}_T(s) = \prod_{p \in \P}  \exp \left( \sum_{n \ge 1} \frac{1}{n} p^{-n s} \right) = \prod_{p \in \P}  \frac{1}{1 - p^{-s}} = \zeta(s) , \]
as we already saw in a slightly different way in Theorem \ref{thm:Riemann_zeta}.
\end{example}

\section{Conclusions}

There are various directions for further study here.    For example, knowing that any tame species $F$ has a Dirichlet series $\widehat{F}$, we can study how various monoidal structures on tame species give operations on Dirichlet series.   We have mentioned four monoidal structures, but another is the \define{substitution product}, where an $F \circ G$ structure on a finite set $S$ is a way of partitioning $S$ into nonempty blocks and then putting a $G$-structure on each block and an $F$-structure on the set of blocks \cite[Sec.\ 1.4]{BLL}.    Is there a multiplicative analogue of this where we write $S$ as a cartesian product of sets of cardinality $> 1$, put a $G$-structure on each factor and an $F$-structure on the set of factors?   What operation on Dirichlet series does this give?

Also, Ramachandran has found fascinating relations between arithmetic zeta functions and the big Witt ring of the integers \cite{Ramachandran}.    Some of his results seem ripe for categorification, especially since we now know how a way to categorify the big Witt comonad---the comonad which when applied to a commutative ring gives its big Witt ring \cite{BMT}.  This categorification is closely related to the substitution product on \textsl{linear} species.   

\section{Acknowledgements}

This paper is an elaboration of joint work done with James Dolan.  We thank Matthew Emerton and other denizens of the $n$-Category Caf\'e for help with zeta functions.

\end{document}